# Two Iterative algorithms for the matrix sign function based on the adaptive filtering technology

Feng Wu[*][1], Keqi Ye [2], Li Zhu [2], Yueling Zhao [3], Jiqiang Hu [3] and Wanxie Zhong [3]

[1] School of Mathematical Sciences, State Key Laboratory of Structural Analysis of Industrial Equipment, Faculty of Vehicle Engineering and Mechanics, Dalian University of Technology, Dalian, 116023, People's Republic of China.

[2] Department of Mechanics, Dalian University of Technology, Dalian 116023, People's Republic of China.

[3] Faculty of Vehicle Engineering and Mechanics, State Key Laboratory of Structural Analysis of Industrial Equipment, Dalian University of Technology, Dalian 116023, People's Republic of China

*Corresponding author:

Tel: 8613940846142; E -mail address: wufeng_chn@163.com

# Abstract

In this paper, two new efficient algorithms for calculating the sign function of the large-scale sparse matrix are proposed by combining filtering algorithm with Newton method and Newton Schultz method respectively. Through the theoretical analysis of the error diffusion in the iterative process, we designed an adaptive filtering threshold, which can ensure that the filtering has little impact on the iterative process and the calculation result. Numerical experiments are consistent with our theoretical analysis, which shows that the computational efficiency of our method is much better than that of Newton method and Newton Schultz method, and the computational error is of the same order of magnitude as that of the two methods.

# 1 Introduction

In the scalar domain, the well-known sign function can be defined as

$$\mathrm{sign}(x) = \begin{cases} 1, & \mathrm{Re}(x) > 0 \\ -1, & \mathrm{Re}(x) < 0 \end{cases}.$$

Roberts extended the sign function from the scalar variable to the matrix variable named as matrix sign function [1], which is an important matrix function with various applications, such as Matrix Lyapunov Equations [2], Algebraic Riccati Equation (ARE) [3-7], Stochastic Differential Equations (SDEs) [4], Yang-Baxter-like equation [8] and so on [9-12].

For the calculation of the matrix sign function, many scholars have carried out related research. The calculation methods for the matrix sign function include direct methods and iterative methods. One representative method of the direct method is the Schur Method which has excellent numerical stability. Let $A \in \mathbb{C}^{n \times n}$ have the Schur decomposition $A = QTQ^*$, where $Q$ is unitary and $T$ is the upper triangular [13], then the Schur Method uses $\mathrm{sign}(A) = Q\mathrm{sign}(T)Q^*$ to evaluate the sign matrix.

Direct methods usually are costly in terms of time and space. Compared with the direct methods, iterative methods need less computing time and storage space, and received a lot of attention. There are many iterative methods for calculating the sign function, among which the more famous ones are Newton's Method, Newton–Schulz Method, Pad´e Family of Iterations, and so on.

In addition, many scholars had constructed some effective new iterative methods recently. As shown in Ref. [14], an iterative method with fourth-order convergence was proposed, and its performance in numerical experiments was discussed. Ref. [4] introduced a new iteration method to calculate the sign function of the matrix, and proved theoretically that the algorithm is stable, and has global convergence. Ref. [8] derived a new iteration scheme for the sign function with fourth-order convergence and introduced the application of sign function in Yang–Baxter-like equation. Ref. [15] proposed a new scheme, which has global convergence with eighth order of convergence. Ref. [16] proposed a dynamics-based approach to develop efficient matrix iterations for computing the matrix sign function. Ref. [17] proposed a direct generalized Steffensen's method for solving the sign function. Ref. [6] present several parallel implementations of the matrix sign function and applied them to the continuous time algebraic Riccati equation.

It should be noted that the above algorithms are still limited to the calculation of the matrix with low dimension at present. Few people have studied the calculation of large-scale matrices. Large sparse matrices will become dense when performing matrix multiplication or matrix inversion, which makes the iterative algorithm need huge storage space and computing time. Although matrix multiplication and inversion will make the matrix dense, there are many elements close to zero, which have little impact on the final calculation results, but will significantly increase the calculation time and storage space.

One way to improve the calculation efficiency is using the filtering algorithm [18,19,21,22]. The filtering algorithm is to construct a threshold and eliminate the elements in the matrix whose absolute value is less than the threshold. Ref. [20] first combined the filtering algorithm with the Newton method to calculate the approximate inverse of the matrix. The filtering algorithm is not only applied to the calculation of matrix inverse, but also extended to the calculation of other matrix functions (such as matrix exponential, matrix inverse and so on) [18,19,21,22].

In this paper, the filtering algorithm is applied to the Newton Method and the Newton Schultz Method to solve the problem of huge storage and time overhead of the large-scale sparse matrices in calculating the sign functions. And the filtering threshold is given theoretically.

The main structure of this paper is as follows. In Section 2, Newton's method (NM) and Newton–Schulz Method (NS) are introduced. In Section 3, according to the error propagation theory, the norm upper bound to the filter matrix and the filtering threshold are given in theory. In Section 4, numerical experiments are carried out on NMF and NSF, and compared with NM and NS.

## 2 Two classical iterative methods

In order to provide a theoretical basis for the filtering algorithm in Section 3, this section introduces some properties of the two classical iterative algorithms. The most widely used and best-known method is the Newton's Method (NM):

$$X_{k+1} = \frac{1}{2}\left(X_k + X_k^{-1}\right), \quad X_0 = A, \tag{1.1}$$

which can be derived in terms of $X^2 = I$. The NM is globally converged and has quadratic convergence. One way to try to remove the inverse in the formula (1.1) is to approximate it by one step of Newton's method for the matrix inverse and we get the Newton–Schulz Method (NS):

$$X_{k+1} = \frac{1}{2} X_k \left(3I - X_k^2\right), \quad X_0 = A,$$

which is multiplication-rich and retains the quadratic convergence of Newton's method. However, it is locally convergent. Only when $\left\|I - A^2\right\| < 1$, the iterative convergence can be guaranteed [13]. Next, a residual recurrence formula is introduced.

**Lemma 2.1.1** The residual terms $R_k = I - X_k^2$ of NM and NS satisfies

$$R_{k+1} = \frac{1}{4} R_k + \frac{1}{4} I - \frac{1}{4}(I - R_k)^{-1}, k = 1, 2, \cdots, \tag{1.2}$$

and

$$R_{k+1} = \frac{3}{4} R_k^2 + \frac{1}{4} R_k^3, k = 1, 2, \cdots, \tag{1.3}$$

respectively. If $\rho(R_k) < 1$, Eq.(1.2) can be rewritten as $R_{k+1} = -\frac{1}{4} R_k^2 (I - R_k)^{-1}$.

**Proof.** Substituting Eq. (1.1) into $R_{k+1} = I - X_{k+1}^2$:

$$\begin{aligned} R_{k+1} &= I - \frac{1}{4} X_k^2 - \frac{1}{2} I - \frac{1}{4} X_k^{-2} \\ &= \frac{1}{2} R_k + \frac{1}{4} X_k^2 - \frac{1}{4} X_k^{-2} \\ &= \frac{1}{2} R_k + \frac{1}{4} I - \frac{1}{4} R_k - \frac{1}{4} X_k^{-2} \\ &= \frac{1}{4} R_k + \frac{1}{4} I - \frac{1}{4}(I - R_k)^{-1}. \end{aligned} \tag{1.4}$$

If $\rho(R_k) < 1$, substituting $(I - R_k)^{-1} = (I + R_k + R_k^2 + \cdots)$ into Eq.(1.4) yields

$$\begin{aligned} R_{k+1} &= \frac{1}{4} R_k + \frac{1}{4} I - \frac{1}{4}(I - R_k)^{-1} \\ &= \frac{1}{4} R_k + \frac{1}{4} I - \frac{1}{4}(I + R_k + R_k^2 + \ldots) \\ &= -\frac{1}{4} R_k^2 (I - R_k)^{-1}. \end{aligned} \tag{1.5}$$

The residual term $R_k = I - X_k^2$ of NS has the following relationship:

$$R_{k+1} = I - X_{k+1}^2 = I - \left[\frac{1}{2}X_k(3I - X_k^2)\right]^2$$
$$= R_k - \frac{5}{4}(I - R_k) + \frac{6}{4}(I - R_k)^2 - \frac{1}{4}(I - R_k)^3 \qquad (1.6)$$
$$= \frac{3}{4}R_k^2 + \frac{1}{4}R_k^3$$

∎

Take the norm of Eq.(1.5):
$$\|R_{n+1}\| \le \frac{1}{4}\|R_n\|^2 \frac{1}{\|I - R_n\|} \qquad (1.7)$$

Take the norm of Eq.(1.6):
$$\|R_{k+1}\| \le \frac{3}{4}\|R_k\|^2 + \frac{1}{4}\|R_k\|^3 \qquad (1.8)$$

# 3. Iterative method with filtering

NM needs one matrix inversion in each iteration, and NS needs two matrix multiplications in each iteration. These two operations will make the sparse matrix dense. Therefore, a very natural idea is to apply the filtering idea to keep the matrix sparse, so as to reduce the amount of calculation and speed up the calculation.

## 3.1 Newton's method with filter

The Newton's method with filter (NMF) can be expressed as Algorithm 3.1.1.

**Algorithm 3.1.1 Newton Iterative algorithm with filtering**

**Input:** $A$ (Initial matrix), $\varepsilon_{tol}$ (Convergence residual)

**Output:** $\text{sign}(A) = \bar{X}_k$ (result)

**begin**

1. $k = 1, \hat{X}_1 = A, e_1 = \|I - \hat{X}_1^2\|$

2. while $e_k > \varepsilon_{tol}$

3. $\bar{X}_{k+1} = \frac{1}{2}(\bar{X}_k + \bar{X}_k^{-1})$

4. $\hat{X}_{k+1} = \bar{X}_{k+1} - F_{k+1}$

5. $\hat{R}_{k+1} = I - \hat{X}_{k+1}^2$; $e_{k+1} = \|\hat{R}_{k+1}\|$

6. $k = k+1$;

7. end while

The fourth step in Algorithm 3.1.1 represents filtering, which means to eliminate the elements in matrix $\bar{X}_k$ which are less than the filtering threshold $\varepsilon_k$ and replace them with zeros, and the $F_k$ represents all filtered elements. Ref. [22] has developed a filtering algorithm to calculate the filtering threshold when the upper bound to $\|F_k\|$ is given. This paper will use this algorithm for the filtering.

Intuitively, the larger $\varepsilon_k$ is, the more elements will be filtered out, which will lead to larger error and even non convergence of calculation. The smaller $\varepsilon_k$ is, the smaller the impact on the result will be, but the calculation time will not be greatly improved. Next, we will discuss how to select the filtering threshold $\varepsilon_k$.

Let us first introduce the following Lemma.

**Lemma 3.1.1** In NMF, if $\|\hat{R}_i\| \ll 1$, $\|F_i\| \ll 1$, and $j \geq i$, then

$$\hat{R}_{j+1} = -\frac{1}{4}\hat{R}_j^2 + \frac{1}{4}\left[F_j \hat{X}_j + F_j \hat{X}_j^{-1} + \hat{X}_j F_j + \hat{X}_j^{-1} F_j\right]. \qquad (1.9)$$

**Proof.** The residual $\hat{R}_{i+1}$ can be written as

$$\hat{R}_{i+1} = I - \hat{X}_{i+1}^2 = I - \frac{1}{4}\left(\hat{X}_i + \hat{X}_i^{-1} - F_i\right)^2.$$

Noting that $\|F_i\| \ll 1$, $F_i^2$ can be ignored, so we have

$$\begin{aligned}
\hat{R}_{i+1} &= I - \frac{1}{4}\left[\hat{X}_i \hat{X}_i + 2I + \hat{X}_i^{-1}\hat{X}_i^{-1} - F_i \hat{X}_i - F_i \hat{X}_i^{-1} - \hat{X}_i F_i - \hat{X}_i^{-1} F_i\right] \\
&= \frac{1}{4}\hat{R}_i + \frac{1}{4}I - \frac{1}{4}\hat{X}_i^{-2} + \frac{1}{4}\left[F_i \hat{X}_i + F_i \hat{X}_i^{-1} + \hat{X}_i F_i + \hat{X}_i^{-1} F\right] \\
&= \frac{1}{4}\hat{R}_i + \frac{1}{4}I - \frac{1}{4}\left(I - \hat{R}_i\right)^{-1} + \frac{1}{4}\left[F_i \hat{X}_i + F_i \hat{X}_i^{-1} + \hat{X}_i F_i + \hat{X}_i^{-1} F\right]
\end{aligned} \qquad (1.10)$$

Because of $\|\hat{R}_i\| \ll 1$, we have

$$\left(I - \hat{R}_i\right)^{-1} = I + \hat{R}_i + \hat{R}_i^2. \qquad (1.11)$$

Substituting Eq. (1.11) into Eq.(1.10) and neglecting higher-order terms yields

$$\hat{R}_{i+1} = -\frac{1}{4}\hat{R}_i^2 + \frac{1}{4}\left[F_i\hat{X}_i + F_i\hat{X}_i^{-1} + \hat{X}_i F_i + \hat{X}_i^{-1} F\right]. \tag{1.12}$$

Replacing the subscript $i$ in Eq. (1.10) with $j$ yields Eq. (1.7). ∎

On the selection of the filtering threshold, a very natural idea is to control the difference between the filtering and non-filtering algorithms, which leads to the adaptive filtering threshold. The residuals involved in the NM and the NMF are $R_k$ and $\hat{R}_k$, respectively. Define the relative residual as $\delta R_k = \hat{R}_k - R_k$.

**Theorem 3.1.2** For the NMF, if $\|\hat{R}_i\| \ll 1$, and

$$\|F_k\| \leq \frac{\|R_k\|^2}{\|\hat{X}_k\| + \|\hat{X}_k^{-1}\|}, k \geq i, \tag{1.13}$$

then $\|\delta R_k\| \leq 4\|R_k\|$.

**Proof.** As $\|\hat{R}_i\| \ll 1$, $\|R_i\|^p$ and $\|\hat{R}_i\|^p$ can be ignored for $p \geq 2$ in the following proof. Noting that $\delta R_1 = \hat{R}_1 - R_1 = 0$ and the Lemma 3.1.1, we have

$$\begin{aligned}\delta R_2 &= \hat{R}_2 - R_2 \\ &= -\frac{1}{4}\left(\hat{R}_1^2 + \hat{R}_1^3 + \ldots\right) + \frac{1}{4}\left(F_1\hat{X}_1 + F_1\hat{X}_1^{-1} + \hat{X}_1 F_1 + \hat{X}_1^{-1} F_1\right) - R_2 \\ &= \frac{1}{4}\left(F_1\hat{X}_1 + F_1\hat{X}_1^{-1} + \hat{X}_1 F_1 + \hat{X}_1^{-1} F_1\right)\end{aligned}$$

By the same way, we can get the recurrence formula of $\delta R_k$:

$$\delta \hat{R}_k = \frac{1}{4}\left(\delta R_{k-1} R_{k-1} + R_{k-1}\delta R_{k-1}\right) + \frac{1}{4}\left(F_{k-1}\hat{X}_{k-1} + F_{k-1}\hat{X}_{k-1}^{-1} + \hat{X}_{k-1} F_{k-1} + \hat{X}_{k-1}^{-1} F_{k-1}\right)$$
$$\tag{1.14}$$

Taking the norm of Eq.(1.14) yields

$$\|\delta R_k\| \leq \frac{1}{2}\|\delta R_{k-1}\|\|R_{k-1}\| + \frac{1}{2}\|F_{k-1}\|\|\hat{X}_{k-1}\| + \frac{1}{2}\|\hat{X}_{k-1}^{-1}\|\|F_{k-1}\|. \tag{1.15}$$

Combining (1.15), (1.7) and (1.13) results in $\|\delta R_k\| \leq 4\|R_k\|$. ∎

Theorem 3.1.2 actually provides the upper bound to $\|F_k\|$.

In the first few steps of the iteration, $\|\hat{R}_i\| \ll 1$ is not satisfied. At this time, a fixed filtering threshold will be recommended, that is $\|F_k\| \leq c_2 \varepsilon_{\text{tol}}$, where $c_2 = 10^{-4}$. Next, we

will analyze what happens to the filtering error $F_k$ in step $k+1$:

$$\begin{aligned}\bar{X}_{k+1} &= \frac{1}{2}(\bar{X}_k + \bar{X}_k^{-1}) = \frac{1}{2}\left[(\hat{X}_k + F_k) + (\hat{X}_k + F_k)^{-1}\right] \\ &= \frac{1}{2}\left[(\hat{X}_k + \hat{X}_k^{-1}) + F_k - \hat{X}_k^{-1} F_k \hat{X}_k^{-1}\right]\end{aligned} \qquad (1.16)$$

At the beginning of the iteration, $\|F_k\| \leq c_2 \varepsilon_{tol}$ is a small quantity, so Eq.(1.16) can be rewritten as $\bar{X}_{k+1} = \frac{1}{2}(\hat{X}_k + \hat{X}_k^{-1})$ which means filtering has little impact to the iteration. When the iteration goes to the later stage, because $\|\hat{R}_i\| \ll 1$, so we can assume $X_k \approx X_k^{-1} \approx \mathrm{sign}(A) = S$, Eq.(1.16) can be rewritten as

$$\bar{X}_{k+1} = \frac{1}{2}\left[(\hat{X}_k + \hat{X}_k^{-1}) + F_k R_k\right]. \qquad (1.17)$$

Because $\|\hat{R}_i\| \ll 1$, then $F_k R_k$ is quite small compared with $\hat{X}_k + \hat{X}_k^{-1}$, so it can be ignored which means NMF will get the same result with NM.

## 3.2 Newton–Schulz Method with filtering

The Newton–Schulz Method with filtering (NSF) can be summarized as Algorithm 3.2.1.

**Algorithm 3.2.1** Newton Schulz iteration with filtering

**Input:** $A$ (Initial matrix), $\varepsilon_{tol}$ (Convergence residual)

**Output:** $\mathrm{sign}(A) = \hat{X}_k$ (result)

**begin**

1. $k = 1, \bar{X}_1 = \hat{X}_1 = A, e_1 = \|I - X_1^2\|$

2. while $e_k > \varepsilon_{tol}$

3. $\bar{X}_{k+1} = \frac{1}{2}\hat{X}_k(3I - \hat{X}_k^2)$

4. $\hat{X}_{k+1} = \bar{X}_{k+1} - F_{k+1}$

4. $\hat{R}_{k+1} = I - \hat{X}_{k+1}^2; e_{k+1} = \|\hat{R}_{k+1}\|$

5. $k = k+1;$

6. end while

**Theorem 3.2.1** when $\|\hat{R}_i\| \ll 1$ for the NSF, if

$$\|F_k\| \leq \frac{\frac{3}{4}\|R_{k-1}\|^2}{3\|\bar{X}_{k-1}\| + \|\bar{X}_{k-1}\|^3}, k \geq i, \tag{1.18}$$

then $\|\delta R_k\| \leq 3\|R_k\|$.

**Proof.** Investigate the change of the Newton Schultz iterative residual with filtering:

$$\begin{aligned}
\bar{R}_{k+1} &= I - \bar{X}_{k+1}^2 \\
&= I - \left(\hat{X}_{k+1} - F_{k+1}\right)^2 \\
&= I - \left[\frac{1}{2}\bar{X}_k\left(3I - \bar{X}_k^2\right) - F_{k+1}\right]^2 \\
&= \frac{\bar{R}_k^3}{4} + \frac{3\bar{R}_k^2}{4} + \frac{3}{2}F_{k+1}\bar{X}_k - \frac{1}{2}F_{k+1}\bar{X}_k^3 + \frac{3}{2}\bar{X}_k F_{k+1} - \frac{1}{2}\bar{X}_k^3 F_{k+1}
\end{aligned} \tag{1.19}$$

Define the relative residuals $\delta R_k = \hat{R}_k - R_k$, noting the Lemma 2.1.1 we have

$$\begin{aligned}
\delta R_k &= \bar{R}_k - R_k \\
&= \frac{\bar{R}_{k-1}^3}{4} + \frac{3\bar{R}_{k-1}^2}{4} + \frac{3}{2}F_k\bar{X}_{k-1} - \frac{1}{2}F_k\bar{X}_{k-1}^3 + \frac{3}{2}\bar{X}_{k-1}F_k - \frac{1}{2}\bar{X}_{k-1}^3 F_k - \frac{3}{4}R_{k-1}^2 - \frac{1}{4}R_{k-1}^3 \\
&= \frac{1}{4}\left[\delta R_{k-1}R_{k-1}R_{k-1} + R_{k-1}\delta R_{k-1}R_{k-1} + R_{k-1}R_{k-1}\delta R_{k-1}\right] \\
&\quad + \frac{3}{4}\left[\delta R_{k-1}R_{k-1} + R_{k-1}\delta R_{k-1}\right] + \frac{3}{2}F_k\bar{X}_{k-1} - \frac{1}{2}F_k\bar{X}_{k-1}^3 + \frac{3}{2}\bar{X}_{k-1}F_k - \frac{1}{2}\bar{X}_{k-1}^3 F_k \\
&= \frac{3}{4}\left[\delta R_{k-1}R_{k-1} + R_{k-1}\delta R_{k-1}\right] + \frac{3}{2}F_k\bar{X}_{k-1} - \frac{1}{2}F_k\bar{X}_{k-1}^3 + \frac{3}{2}\bar{X}_{k-1}F_k - \frac{1}{2}\bar{X}_{k-1}^3 F_k
\end{aligned}$$

$$\tag{1.20}$$

Take the norm of Eq.(1.20),

$$\|\delta R_k\| \leq \frac{3}{2}\|\delta R_{k-1}\|\|R_{k-1}\| + 3\|F_k\|\|\bar{X}_{k-1}\| + \|F_k\|\|\bar{X}_{k-1}\|^3 \tag{1.21}$$

Combining (1.21), (1.8) and (1.18) results in $\|\delta R_k\| \leq 3\|R_k\|$. ∎

Filtering is equivalent to artificially introducing some errors. An important consideration is how the artificially introduced errors will propagate with the iterative process. Next, we will analyze what happens to the filtering error $F_k$ in step $k+1$

$$\begin{aligned}
\bar{X}_{k+1} &= \frac{1}{2}\hat{X}_k\left(3I - \hat{X}_k^2\right) \\
&= \frac{1}{2}(\bar{X}_k - F_k)\left[3I - (\bar{X}_k - F_k)^2\right] \\
&= \frac{1}{2}\bar{X}_k\left[3I - (\bar{X}_k - F_k)^2\right] - \frac{1}{2}F_k\left[3I - (\bar{X}_k - F_k)^2\right] \\
&= \frac{1}{2}\bar{X}_k\left[3I - \bar{X}_k^2 + 2\bar{X}_k F_k - F_k^2\right] - \frac{1}{2}F_k\left[3I - \bar{X}_k^2 + 2\bar{X}_k F_k - F_k^2\right] \\
&= \frac{1}{2}\bar{X}_k\left[3I - \bar{X}_k^2 + 2\bar{X}_k F_k\right] - \frac{1}{2}F_k\left[3I - \bar{X}_k^2\right] \\
&= \frac{1}{2}\bar{X}_k\left[3I - \bar{X}_k^2\right] + \frac{3}{2}\bar{X}_k\bar{X}_k F_k - \frac{3}{2}F_k \\
&= \frac{1}{2}\bar{X}_k\left[3I - \bar{X}_k^2\right] + \frac{3}{2}(\bar{X}_k\bar{X}_k - I)F_k \\
&= \frac{1}{2}\bar{X}_k\left[3I - \bar{X}_k^2\right] - \frac{3}{2}R_k F_k
\end{aligned} \quad (1.22)$$

Therefore, at the beginning of the iteration $\|R_k\| \sim 1$, we adopt a fixed threshold $\|F_k\| \leq c_4 \varepsilon_{\text{tol}}$, where $c_4 = 10^{-4}$ is a parameter. Because $F_k$ is a matrix with small elements and norm, so the error matrix $\frac{3}{2}F_k R_k$ in the next step can be ignored. In the later stage $\|R_k\| \ll 1$, it can be seen from Eq.(1.22) that the error $F_k$ introduced in step $k$ becomes $\frac{3}{2}F_k R_k$ in step $k+1$, which is small. In other words, the filtering error can be ignored.

# 4 Numerical Experiment

This section will test several algorithms (NMF, NSF) introduced before. In the calculation, we used the MATLAB's sparse matrix storage technology (sparse). When the residual converges to $\|\hat{R}_k\| < 10^{-6}$, the filtering threshold of NMF is determined by Eq.(1.13), and the filtering threshold of NSF is determined by Eq.(1.18), otherwise the fixed filtering threshold $\|F_k\| \leq 10^{-4}\varepsilon_{\text{tol}}$ was used for both NMF and NSF. The CPU used is Intel ® Core ™ I7-10750h processor, and the memory size is 32GB. The program has been uploaded to the website: https://rocewea.com/6.html, and the matrices used in 4.4 can also be found on this website.

For the convenience of the description, hereinafter, the execution times of different algorithms is abbreviated as ET, the absolute error（abbreviated as AE）is defined by

$\|X_f - \hat{X}_f\|$, which represents the difference between the NM and NMF results (or the NS and NSF results).

## 4.1 Filter error test

The example in this subsection is used to investigate the relationship between $\|\hat{X}_k\|$ and $\|X_k\|$ and that between $\|\hat{R}_k\|$ and $\|R_k\|$.

First, we used the matrix randomly generated by MATLAB code $A = \text{sprand}(n, n, 0.01)$, where $n$ is the dimension. Because the generated matrix may not meet the convergence condition of NSF, NS method, we used the NMF and NM methods with $\varepsilon_{\text{tol}} = 10^{-12}$ for this matrix. Three dimensions, i.e., $n = 500$, 1000, and 2000, are used. For each dimension, the computation is performed 5 times. The norm of $\delta X_k$ and $\delta R_k$ for different dimensions are plotted in Fig. 4.4.1.

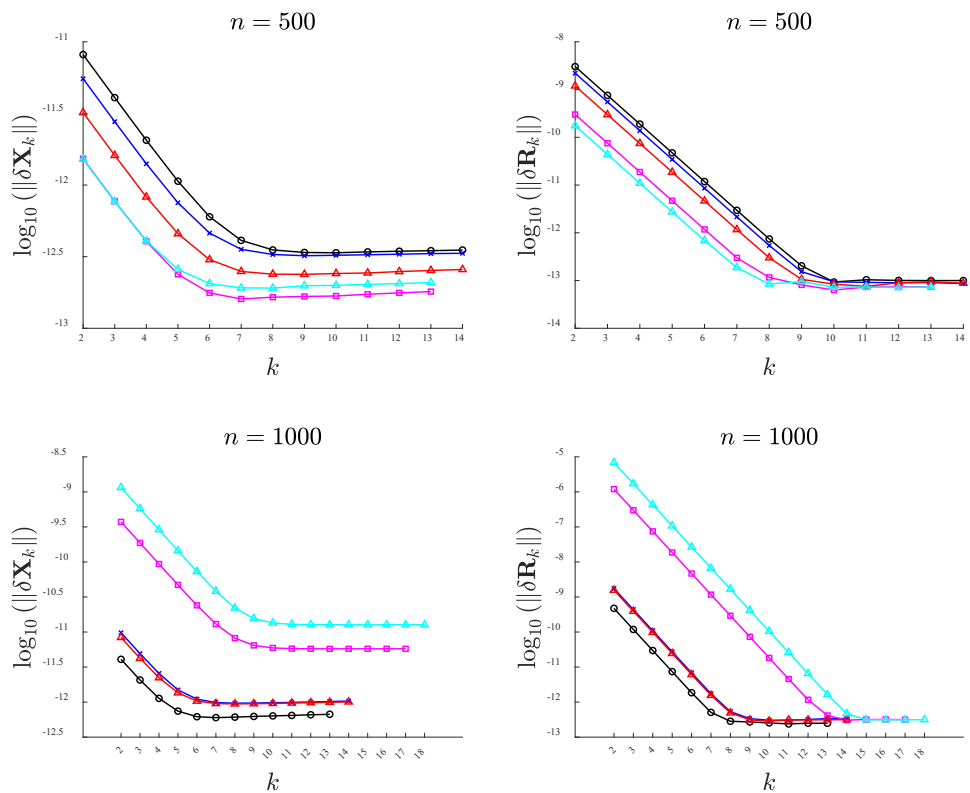

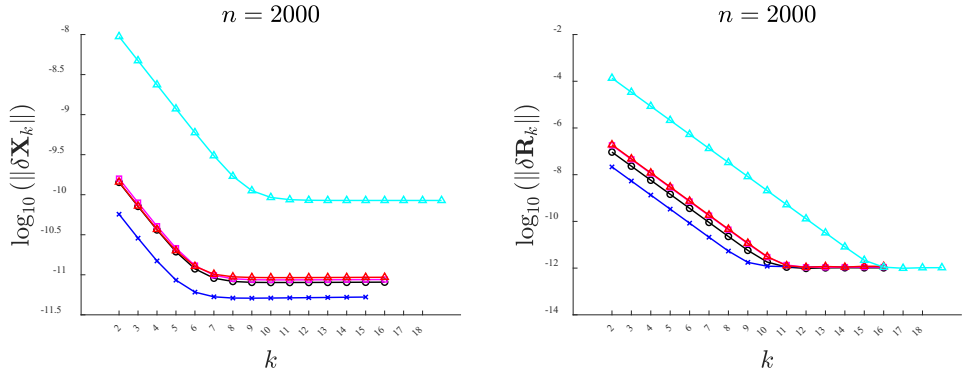

Fig 4.1.1 The norm of $\delta X_k$ and $\delta R_k$ for NM and NMF with different dimensions.

From Fig 4.1.1, it can be seen that, for the NMF, the norm of $\delta X_k$ and $\delta R_k$ becomes smaller and smaller as the iteration goes, which shows the filtering idea applied to NM has little impact on the convergence result. This is very consistent with the discussion on the selection of filtering threshold in section 3.

The matrix used next is the positive definite symmetric matrix randomly generated by MATLAB code: $A = \text{sprandsym}(n, 0.01, 1)$, where $n$ is the dimension of the matrix. The algorithm NSF and NS with $\varepsilon_{\text{tol}} = 10^{-12}$ are used for this matrix. Three dimensions, i.e., $n = 500$, 1000, and 2000, are used. For each dimension, the computation is performed 5 times. The norm of $\delta X_k$ and $\delta R_k$ for different dimensions are plotted in Fig. 4.1.2.

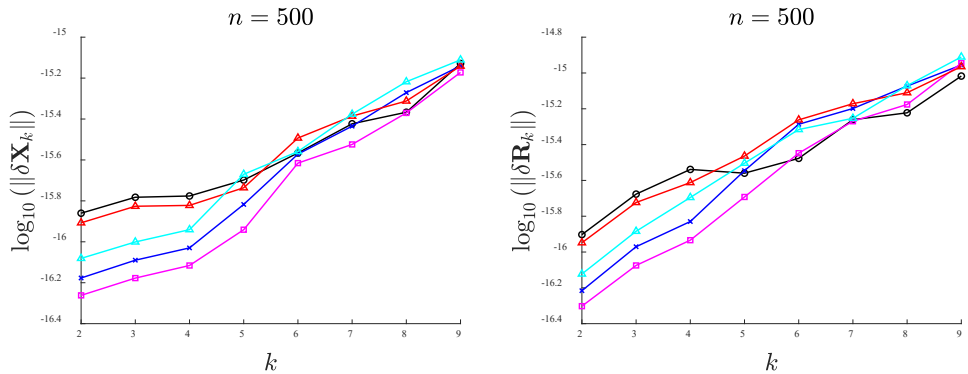

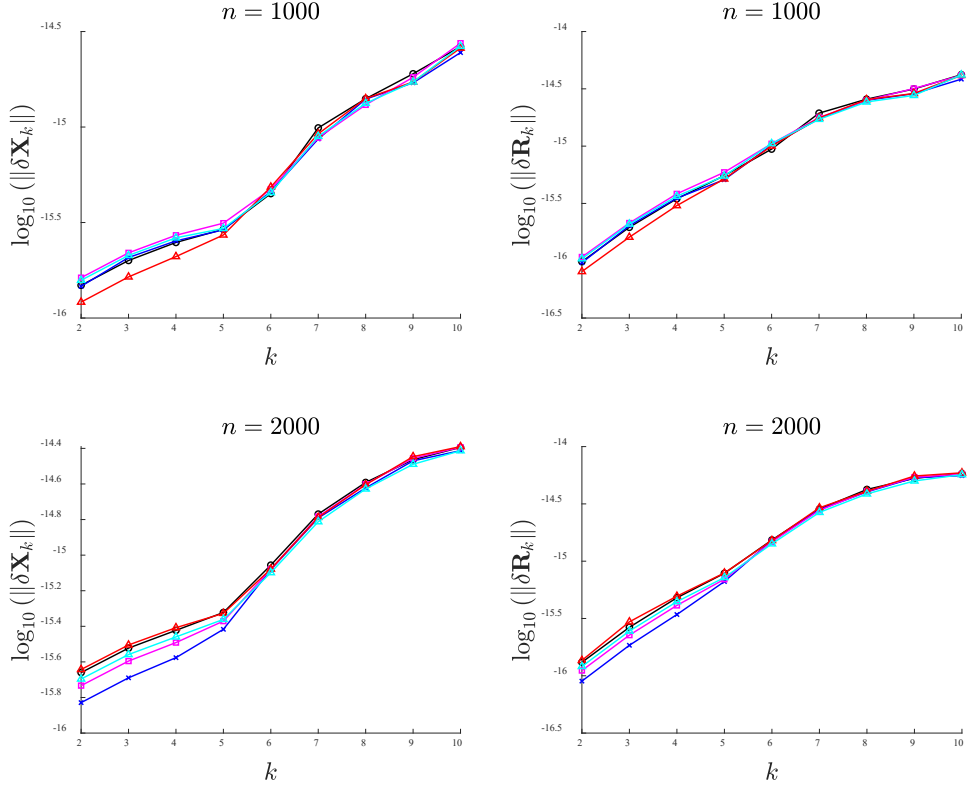

Fig 4.1.2 The norm of $\delta\boldsymbol{X}_k$ and $\delta\boldsymbol{R}_k$ for NS and NSF with different dimensions.

From Fig 4.1.2, it can be seen that, for the NSF, both the norm of $\delta\boldsymbol{X}_k$ and $\delta\boldsymbol{R}_k$ are smaller than the given $\varepsilon_{\text{tol}} = 10^{-12}$ during the iterative process. The numerical result shows that the filtering idea has also little impact on the convergence result of the NM.

## 4.2 Matrix differential equation

In this section, we compared the performance of NMF, NSF, NM and NS with $\varepsilon_{\text{tol}} = 10^{-12}$ in solving the following matrix differential equations

$$\boldsymbol{y}''(t) + \boldsymbol{B}\boldsymbol{y}(t) = 0. \tag{1.23}$$

Equation (1.23) has the initial values $\boldsymbol{y}(0) = \boldsymbol{\alpha}$ and $\boldsymbol{y}'(0) = \boldsymbol{\beta}$. The equation has an analytical solution, which can be expressed as

$$\boldsymbol{y}(t) = \boldsymbol{\alpha}\cos\left(\boldsymbol{B}^{1/2}t\right) + \boldsymbol{\beta}\boldsymbol{B}^{-1/2}\sin\left(\boldsymbol{B}^{1/2}t\right). \tag{1.24}$$

Obviously, to get the solution of the problem, we need to calculate $\boldsymbol{B}^{1/2}$, $\boldsymbol{B}^{-1/2}$ and consider the matrix:

$$B = \text{spdiags}\left(\text{kron}\left(\text{ones}(n,1),[\lambda,1-2\lambda,\lambda]\right),[-1,0,1],n,n\right), \quad (1.25)$$

where $\lambda = \dfrac{1}{16}$. The matrix is related to the finite difference dispersion of the spatial second derivative [9]. According to the Theorem 6.1 in Ref. [4]:

$$A = \begin{bmatrix} 0 & B \\ I & 0 \end{bmatrix}, \quad \text{sign}(A) = \begin{bmatrix} 0 & B^{1/2} \\ B^{-1/2} & 0 \end{bmatrix} \quad (1.26)$$

Therefore, it is very convenient to obtain $B^{1/2}$, $B^{-1/2}$ at one time in terms of (1.26). This section used the NM, NS, NMF, and NSF to calculate $\text{sign}(A)$ with a series of dimensions. ET and AE of these four methods are listed in Table 4.2.1.

Table 4.2.1 ET and AE of different algorithms

| $n$ | NMF | | NM | NSF | | NS |
|---|---|---|---|---|---|---|
| | ET (s) | AE | ET (s) | ET | AE | ET (s) |
| 500 | 0.0589 | 7.57E-14 | 0.2286 | 0.0098 | 6.41E-14 | 0.0189 |
| 1000 | 0.1407 | 1.08E-13 | 2.1514 | 0.011 | 9.21E-14 | 0.0346 |
| 3000 | 0.9482 | 1.88E-13 | 25.1071 | 0.0268 | 1.61E-13 | 0.1275 |
| 5000 | 2.4548 | 2.42E-13 | 64.1395 | 0.0481 | 2.08E-13 | 0.2187 |

It is not difficult to see from Table 4.2.1 that the calculation time of NSF is the shortest, and the calculation time of NMF is greater than that of NS, which is due to the matrix inverse usually costs many times. Due to the addition of the filtering idea, the ET of the NMF and NSF decreases significantly, and the improvement is more and more obvious with the increase of the matrix dimension. In terms of the calculation error, the absolute error of NMF and NSF remains at a very low level.

## 4.3 Solving algebraic Riccati equation

Sign functions also have important applications in solving the algebraic Riccati equation (ARE). ARE is mainly used to solve linear quadratic optimization problems in control theory. For the selected linear quadratic performance index, we can get the corresponding Riccati equation. If the Riccati equation has a positive definite solution, the linear quadratic optimization problem has a solution, and we can get the corresponding state feedback. ARE has the following forms

$$P(U) = UC + C^T U + Q - UBD^{-1}B^T U = 0, \quad (1.27)$$

where $B \in \mathbb{C}^{n \times n}, C \in \mathbb{C}^{n \times m}$, $Q$ is a positive semidefinite matrix of $n$ dimension and $D$ is a positive definite matrix of $n$ dimension. $U$ is the unknown matrix. According to Ref. [4], there is

$$S = \text{sign}(A) = \text{sign}\left(\begin{bmatrix} C & BD^{-1}B^T \\ Q & -C^T \end{bmatrix}\right) = \begin{bmatrix} W_{11} & W_{12} \\ W_{21} & W_{22} \end{bmatrix}; \quad S\begin{bmatrix} I \\ -U \end{bmatrix} = \begin{bmatrix} -I \\ U \end{bmatrix} \quad (1.28)$$

Once $A$ is obtained, $U$ can be computed by

$$U = W_{12}^{-1}(W_{11}+I) \tag{1.29}$$

Therefore, we can solve ARE according to Eq.(1.29) by calculated the sign function matrix according to Eq.(1.28). Let

$$\begin{aligned} B &= \text{spdiags}\left(\text{kron}(\text{ones}(n,1),[-1.6,0,0.8,0,-1.6]),[-2,-1,0,1,2],n,n\right) \\ C &= \text{spdiags}\left(\text{kron}\left(\text{ones}(n,1),\left[\frac{1}{16},\frac{7}{8},\frac{1}{16}\right]\right),[-1,0,1],n,n\right) \end{aligned} \tag{1.30}$$

$D$ is a random symmetric positive definite matrix, and $Q = BD^{-1}B^T$. In this example, inv function of MATLAB is used for inversion. The NMF, NM, NSF, NS methods with $\varepsilon_{tol} = 10^{-14}$ are used for this problem. The signm function in MATLAB is also used for this example. The ET, and equation errors (EE) of different methods are listed in Table 4.3.1. The Equation error (EE) is defined by $\|P(U^*)\|_\infty$, where $U^*$ represents the matrix solution based on different methods.

Table 4.3.1 ET (s) and EE of different algorithms

| n | NMF | | NM | | NSF | | NS | | signm | |
|---|---|---|---|---|---|---|---|---|---|---|
| | ET(s) | EE | ET(s) | EE | ET(s) | EE | ET(s) | EE | ET(s) | EE |
| 500 | 3.3 | 2.0E-06 | 4.2 | 1.9E-06 | 2.6 | 9.6E-06 | 21.3 | 3.0E-06 | 55.2 | 1.2E-04 |
| 600 | 3.9 | 1.4E-06 | 6.1 | 1.6E-06 | 3.3 | 3.2E-06 | 27.2 | 3.1E-06 | 90.7 | 6.3E-04 |
| 700 | 6.6 | 2.7E-05 | 12.5 | 2.8E-05 | 4.7 | 1.3E-05 | 51.5 | 5.2E-06 | 139.4 | 3.1E-04 |

It can be seen from Table 4.3.1 that due to the addition of the filtering idea, the amount of calculation is significantly reduced, and the improvement is more and more obvious with the increase of the matrix dimension. Compared with the signm, the NMF and NSF perform better in terms of the computational time and accuracy.

## 4.4 Calculation of large-scale sparse matrix

The fourth example employ 31 large-scale sparse matrices to test the performance of the NMF and NSF. All these matrices can be divided into three types. Type 1 matrixes are the coefficient matrices involved in the function approximation problems, which can be expressed as $\int_\Omega NN^T dxdy$, in which $N$ is the vector consisting of compactly supported functions. Type 2 are matrices used in the finite element dynamic analysis which can be written as $M + \gamma \Delta t C + \beta \Delta t^2 K$, where $M$, $C$, $K$ is the mass matrix, damping matrix and stiffness matrix, respectively, $\Delta t$ the time step and $\beta, \gamma$ the integral constants. Type 3 are encountered in calculating the communicability and the node centralities of complex networks, which are the form of $I-\alpha H$, where $H$ is the adjacent matrix. All the adjacent matrix are download from the Network Repository at: http://networkrepository.com/.

In this section, the NMF and NSF methods with $\varepsilon_{tol} = 10^{-13}$ will be used to calculate the sign functions of these matrices. Since the dimensions of these matrices are between 10000 and 1508065 dimensions, the inversion function inv of MATLAB takes a long time, so the code inv_filter [23] is used for the inversion in this subsection. We recorded the ET and the residual, defined by $\|\hat{R}_f\|$, where f is the last iteration step.

It can be seen from the Table 4.4.1 that the NMF and NSF compute the sign of these large and sparse matrices efficiently and precisely. For all the 31 matrices, the residuals $\|\hat{R}_f\|$ are smaller than the given $\varepsilon_{tol} = 10^{-13}$.

Table 4.4.1 ET (s) and $\|\hat{R}_f\|$ of NMF and NSF

| No | Type | n | NMF | | NSF | |
|---|---|---|---|---|---|---|
| | | | ET(s) | $\|\hat{R}_f\|$ | ET(s) | $\|\hat{R}_f\|$ |
| 1 | Type 1 | 10000 | 2.048 | 5.85E-17 | 0.735 | 3.72E-15 |
| 2 | Type 1 | 10000 | 5.737 | 2.24E-15 | 1.489 | 2.22E-14 |
| 3 | Type 1 | 10000 | 14.675 | 1.71E-16 | 2.745 | 1.09E-15 |
| 4 | Type 1 | 10000 | 14.675 | 1.71E-16 | 5.712 | 6.80E-14 |
| 5 | Type 1 | 10000 | 6.932 | 7.45E-16 | 1.371 | 7.36E-16 |
| 6 | Type 3 | 12406 | 0.068 | 4.71E-17 | 0.079 | 1.27E-17 |
| 7 | Type 2 | 18150 | 89.908 | 3.89E-16 | 28.456 | 1.47E-14 |
| 8 | Type 3 | 19502 | 1.371 | 4.86E-17 | 0.315 | 1.46E-14 |
| 9 | Type 3 | 25626 | 0.33 | 9.82E-21 | 5.62 | 2.62E-16 |
| 10 | Type 2 | 28212 | 46.12 | 2.65E-16 | 6.883 | 1.93E-14 |
| 11 | Type 3 | 31022 | 1.269 | 5.05E-14 | 0.446 | 3.14E-16 |
| 12 | Type 3 | 31385 | 2.54 | 4.16E-16 | 0.752 | 4.58E-14 |
| 13 | Type 3 | 32075 | 2.54 | 4.16E-16 | 1.359 | 6.86E-14 |
| 14 | Type 3 | 38120 | 0.545 | 1.71E-20 | 12.791 | 2.29E-16 |
| 15 | Type 3 | 38744 | 1.122 | 4.90E-15 | 3.9 | 2.22E-16 |
| 16 | Type 3 | 57735 | 4.823 | 7.48E-15 | 3.637 | 2.22E-16 |
| 17 | Type 3 | 68121 | 1.083 | 4.06E-15 | 0.566 | 2.51E-16 |
| 18 | Type 3 | 73283 | 7.014 | 3.91E-16 | 1.358 | 3.21E-14 |
| 19 | Type 3 | 83995 | 10.948 | 2.55E-16 | 2.572 | 3.22E-14 |
| 20 | Type 3 | 185639 | 1.633 | 6.36E-15 | 77.518 | 2.22E-16 |
| 21 | Type 3 | 203954 | 13.161 | 7.34E-16 | 3.87 | 5.50E-14 |
| 22 | Type 3 | 245874 | 0.212 | 5.20E-16 | 0.217 | 2.20E-13 |
| 23 | Type 3 | 343791 | 0.467 | 4.07E-15 | 0.527 | 2.22E-16 |
| 24 | Type 3 | 381689 | 0.998 | 4.38E-15 | 8.844 | 2.22E-16 |
| 25 | Type 3 | 503625 | 4.659 | 3.55E-15 | 21.684 | 2.23E-16 |
| 26 | Type 3 | 503625 | 4.651 | 5.09E-15 | 19.612 | 2.23E-16 |
| 27 | Type 3 | 504855 | 1.912 | 1.74E-17 | 13.776 | 2.61E-16 |
| 28 | Type 3 | 504855 | 1.981 | 1.93E-17 | 13.537 | 2.60E-16 |
| 29 | Type 3 | 504855 | 0.728 | 6.65E-19 | 3.425 | 2.24E-16 |
| 30 | Type 3 | 943695 | 1.642 | 6.28E-15 | 3.724 | 5.05E-14 |
| 31 | Type 3 | 1508065 | 1.93 | 3.66E-15 | 5.45 | 2.22E-16 |

# 5 Conclusion

In this paper, the filtering algorithm is combined with Newton method and Newton Schultz method, and two new calculation methods, denoted by NMF and NSF, are proposed. Theoretical analysis shows that the error caused by each filtration on the NMF and NSF has little impact to the final result.

Examples 1 show that the filtering algorithm has little impact on the convergence result of the NMF and NSF. Examples 2 and 3 tested the performance of NMF and NSF in the practical application of sign function. Compared with the original algorithm (NM, NS), NMF and NSF has great efficiency advantages and the error control is still at a very low level. The last example proves that NSF and NMF can compute efficiently and accurately the sign functions of large-scale sparse matrices which are nearly sparse.